\documentclass[reqno]{amsart}
\usepackage{amsmath,amssymb,cite}
\usepackage{graphics,epsfig,color}
\usepackage{orcidlink}
\usepackage{geometry,bm,enumitem}
\usepackage[mathscr]{euscript}
\usepackage[english]{babel}
\usepackage{tikz}
\usepackage{pgfplots}
\usepackage{subcaption} 
\pgfplotsset{compat=1.18}

\theoremstyle{plain}
\newtheorem{theorem}{Theorem}[section]

\theoremstyle{remark}
\newtheorem{remark}[theorem]{Remark}

\numberwithin{equation}{section}
\numberwithin{theorem}{section}

\newcommand{\mc}[1]{{\mathcal #1}}
\newcommand{\bb}[1]{{\mathbb #1}}
\newcommand{\rmd}{\mathrm{d}}
\newcommand{\eps}{\varepsilon}

\title[From kinetic activation dynamics to Holling type predator--prey systems]{From kinetic activation dynamics to Holling type reaction--diffusion predator--prey systems}

\author[Paolo Butt\`a]{Paolo Butt\`a \orcidlink{0000-0002-3193-4282}}
\address{Paolo Butt\`a\hfill\break \indent
Dipartimento di Matematica, Sapienza Universit\`a di Roma 
\hfill\break \indent
P.le Aldo Moro 5, 00185 Roma, Italy}
\email{butta@mat.uniroma1.it}

\date{\today}

\begin{document}

\begin{abstract}
We consider a kinetic model of predator--prey dynamics in which predation is characterized by a sequence of behavioral transitions. Prey and predators are described by distribution functions in phase space, with the predator population further structured into internal states corresponding to inactive, active, and saturated individuals. We show that different macroscopic functional responses emerge from distinct mechanisms governing predator activation. More precisely, activation driven uniquely by pairwise encounters leads to a saturating Holling type~II response, whereas activation additionally modulated by local prey density produces a generalized sigmoidal Holling type III response. This is achieved by means of a diffusive scaling and a formal hydrodynamic limit, from which we derive a closed reaction--diffusion system of Rosenzweig--MacArthur type for the prey and total predator spatial densities.
\end{abstract}



\maketitle

\section{Introduction}
\label{sec:1}

The model developed in this work is based on the idea that predation is not an instantaneous event but unfolds through a sequence of behavioral stages. This viewpoint goes back to the seminal contribution of Holling \cite{Holling1959}, who introduced the concept of handling time and formulated the functional responses that have become central in theoretical ecology. In particular, Holling’s type II response underlies the classical Rosenzweig--MacArthur predator--prey model \cite{RosenzweigMacArthur1963}. Subsequently, numerous models have been proposed to explain the emergence of these functional responses under different ecological conditions.

Among them, individual-based, stochastic, and mean-field approaches have been used to recover Holling’s type II and III responses \cite{DawesSouza2013}. Similarly, approaches based on structured reaction--diffusion systems have been proposed to obtain saturating nonlinearities at the macroscopic scale \cite{ConfortoDesvillettesSoresina2017}. In contrast, we formulate a kinetic model that explicitly integrates spatial transport, velocity dynamics, and predator internal states.

Kinetic approaches to interacting populations have been developed to connect microscopic interaction rules with macroscopic population dynamics, as shown for instance in kinetic predator--prey models leading to Lotka--Volterra-type systems through suitable scaling limits \cite{BondesanZanella2025}. However, the role of internal behavioral activation mechanisms in determining emergent functional responses has received comparatively less attention within such frameworks.

From an individual-level perspective, empirical studies indicate that predators typically pass through intermediate behavioral states before reaching full hunting efficiency. Different processes, such as prey recognition, learning, attention allocation, and adaptation to prey defenses, can contribute to these transitions. These mechanisms are naturally described within structured population frameworks, where behavioral or physiological traits evolve alongside spatial movement. Kinetic transport equations are particularly suitable to this setting, as they allow the simultaneous representation of transport in physical space and progression through internal states. For a general introduction to the kinetic theory of biological systems, we refer to Perthame’s monograph \cite{Perthame2007}.

Motivated by these considerations, we describe the predator population using three internal states: an inactive state, in which individuals are present but unable to capture prey; an active state, in which predators can successfully attack; and a saturated state, representing handling or digestion following a capture. Transitions between predator states are driven by prey encounters and intrinsic relaxation mechanisms.

Compartmental descriptions of predation processes have a long tradition at the population level. In particular, early mechanistic models show how multi-stage dynamics, with distinct searching, attacking, and handling phases, can give rise to Holling type functional responses~\cite{Real1977}. More generally, nonlinear functional responses are often interpreted as emergent effects of hidden behavioral states and internal transitions occurring at the individual level, see, e.g.,~\cite{Berardo2020}. Related two-state formulations, such as searching/handling predator models, also appear in the framework of \cite{Beddington1975,DeAngelis1975}, as well as in stochastic and individual-based approaches to predator--prey interactions.

While these approaches are mainly formulated in terms of ordinary differential equations or discrete stochastic processes, we incorporate these behavioral transitions into a kinetic description and show that two different activation mechanisms give rise to different macroscopic functional responses. To this end, we consider three characteristic timescales: a fast timescale associated with velocity relaxation, an intermediate timescale governing transitions between predator internal states, and a slow timescale associated with the evolution of local population densities. Under a diffusive scaling, we formally derive a closed reaction--diffusion system of Rosenzweig--MacArthur type governing the macroscopic evolution of the prey and total predator spatial densities. The resulting functional response is of Holling type~II when activation is solely driven by pairwise encounters, whereas a density-dependent activation rate leads to a generalized type~III response.

Reaction--diffusion predator--prey systems play a central role in spatial ecology and mathematical biology \cite{Murray1989}; in particular, the class of models obtained here has generated a rich literature, including the analysis of traveling-wave phenomena (see, e.g., \cite{DucrotLiuMagal2021,Dunbar1986,LiWu2008}).

The paper is organized as follows. In Section \ref{sec:2}, we introduce the kinetic model for the interacting species, detailing the internal activation states and the continuous velocity dynamics driven by the Fokker-Planck operators. Section \ref{sec:3} is devoted to the hydrodynamic limit and the internal-state reduction; here, we introduce the diffusive scaling, perform a formal Hilbert expansion, and derive the macroscopic target systems featuring the emergent Holling type functional responses. Section \ref{sec:4} provides a discussion of the results and some concluding remarks.

\section{Kinetic formulation}
\label{sec:2}

We now introduce the kinetic model that encodes the individual-level mechanisms described in the Introduction. The prey population is described by the distribution function $f(x,v,t)$, where $x\in\mathbb R^d$ and $v\in\mathbb R^d$ denote position and velocity, respectively, and where $d=1,2,3$ is the dimension of the physical space. The predators are instead described by three distribution functions $g_i(x,v,t)$, $i=0,1,2$, associated with the inactive, active, and saturated states, respectively. All distributions are nonnegative and integrable with respect to $(x,v)$. In particular, the phase-space integral of $f$ represents the total population size of prey, while the phase-space integrals of $g_i$ represent the total population sizes of predators in the corresponding internal states.

The kinetics dynamics prescribes transport in physical space, stochastic fluctuations of velocity, and transitions between internal predator states triggered by encounters with prey. Movement is modeled through a Fokker--Planck diffusion in velocity, which drives rapid relaxation toward local Maxwellians and provides the natural fast timescale for hydrodynamic closure.

To model predator behavioral transitions, we introduce two alternative activation mechanisms, reflecting distinct ecological interpretations. The first assumes a constant baseline activation rate $\alpha$, representing opportunistic predators that become activated following a single encounter with prey and leading, at the macroscopic level, to a Holling type II functional response. The second assumes a density-dependent activation rate, accounting for predator learning, search image formation, or the progressive saturation of prey refuges, and yielding a sigmoidal Holling type III functional response. We refer to these mechanisms as ``opportunistic activation'' and ``learning-driven activation'', respectively. At the microscopic level, this distinction is encoded through the activation term. In the opportunistic case, activation is modeled by the standard bilinear contact term $\alpha f g_0$. In the learning-driven case, activation is described by the hybrid microscopic-macroscopic interaction term $\alpha \rho_f f g_0$, where $\rho_f(x,t) := \int\! \rmd v \, f(x,v,t)$ is the macroscopic (or local) prey density.

Encounters between prey and active predators, represented by a bilinear interaction term $a f g_1$, with $a>0$, produce both prey losses and demographic recruitment of predators. The choice of a bilinear form is consistent with a low-density population regime, in which pairwise encounters provide the leading-order contribution to predator-prey interactions.

We emphasize that the local dependence on both position and velocity in the bilinear contact terms $fg_0$ and $fg_1$ is not meant to represent exact velocity matching between prey and predators. Rather, it should be viewed as an effective encounter term describing interactions between individuals that are spatially close and share similar movement characteristics.

We refer to Remark~\ref{rem:2.2} below for further discussions regarding the phase-space consistency of these choices. The resulting kinetic equations are given in \eqref{eq:prey} and \eqref{eq:pred1}--\eqref{eq:pred2} below.

Prey dynamics are governed by a transport-diffusion equation with logistic growth and losses due to encounters with active predators:
\begin{equation}
\label{eq:prey}
\partial_t f + v\cdot\nabla_x f = \mc L_1 f + r f \Big(1-\frac{\rho_f}{K}\Big) - a f g_1.
\end{equation}
Here, 
\begin{equation}
\label{mcL1}
\mc L_1 f :=  D_1  \Delta_vf + \nabla_v \cdot (vf) 
\end{equation}
is the Fokker-Planck operator with diffusion coefficient $D_1>0$, while $\rho_f(x,t)$ is the aforementioned local prey density, and $a>0$ the predation rate mediating prey mortality through active predators. Note that the logistic term acts uniformly across the velocity distribution, since demographic growth and intra-specific competition depend on the local population density $\rho_f(x,t)$ and not on velocity.

Predator dynamics couple transport, diffusion, and internal transitions, and we consider the two distinct already mentioned behavioral scenarios:
\\ Case I: Opportunistic activation
\begin{equation}
\label{eq:pred1}
\begin{aligned}
\partial_t g_0 + v\cdot\nabla_x g_0 & = \mc L_2  g_0 - \alpha f g_0 + h^{-1} g_2 - m g_0 + e a f g_1, \\ \partial_t g_1 + v\cdot\nabla_x g_1 & = \mc L_2  g_1 + \alpha f g_0 - \beta f g_1 - m g_1, \\ \partial_t g_2 + v\cdot\nabla_x g_2 & = \mc L_2  g_2 + \beta f g_1 - h^{-1} g_2 - m g_2.
\end{aligned}
\end{equation}
\\ Case II: Learning-driven activation
\begin{equation}
\label{eq:pred2}
\begin{aligned}
\partial_t g_0 + v\cdot\nabla_x g_0 & = \mc L_2  g_0 - \alpha\rho_f f g_0 + h^{-1} g_2 - m g_0 + e a f g_1, \\ \partial_t g_1 + v\cdot\nabla_x g_1 & = \mc L_2  g_1 + \alpha \rho_f f g_0 - \beta f g_1 - m g_1, \\ \partial_t g_2 + v\cdot\nabla_x g_2 & = \mc L_2  g_2 + \beta f g_1 - h^{-1} g_2 - m g_2.
\end{aligned}
\end{equation}
Analogously to \eqref{mcL1},
\begin{equation}
\label{mcL2}
\mc L_2 g := D_2 \Delta_vg + \nabla_v \cdot (vg)
\end{equation}
is the Fokker-Planck operator with diffusion coefficient $D_2>0$. In both cases, inactive predators $g_0$ become active at rate $\alpha f$ (case I) or $\alpha \rho_f f$ (case II); active predators $g_1$ capture prey and enter the saturated state $g_2$ at rate $\beta f$. Additionally, predation events contribute to predator recruitment into the inactive class via the term $e a f g_1$, where $e>0$ denotes the conversion efficiency, i.e., the fraction of consumed prey that is converted into predator recruitment. Saturated predators return to the inactive state at rate $h^{-1}$ ($h$ being the ``handling time"). All states are subject to a basal mortality $m>0$.

\begin{remark}[On the choice of the Ornstein--Uhlenbeck operator] \label{rem:2.1} 

We model the unbiased random fluctuations in velocity, arising from behavioral variability in prey and predator movements, using the Fokker--Planck operators $\mc L_1$ and $\mc L_2$ rather than simple velocity diffusions $D_1\Delta_v$ and $D_2\Delta_v$, as the latter lack a mechanism driving the velocity distribution toward an equilibrium profile. 

For simplicity, we also assume that the three predator states share the same velocity diffusion operator. Clearly, one may expect inactive, active, and handling predators to exhibit different mobility patterns. However, allowing state-dependent diffusion coefficients would generally lead, after elimination of the fast behavioral variables, to more complex effective diffusion terms for the total predator density, possibly involving nonlinear or cross-diffusive contributions. The present assumption represents the simplest setting in which the hydrodynamic limit closes into a standard reaction--diffusion system of Rosenzweig--MacArthur type.  
\end{remark}

\begin{remark}[On the hybrid microscopic-macroscopic  activation mechanism]
\label{rem:2.2}

The activation mechanism in Case~II is modeled by the hybrid interaction term $\alpha \rho_f(x,t)\, f(x,v,t)\, g_0(x,v,t)$, which combines a local phase-space interaction with a macroscopic density factor. While this structure departs from purely bilinear kinetic interactions, it is introduced to capture a fundamental feature of predator behavior, namely that activation is not solely triggered by individual encounters, but is modulated by the overall environmental availability of prey.

From a modeling perspective, a fully local cubic interaction term of the form $f^2 g_0$ would represent the simultaneous encounter of an inactive predator with two prey sharing the same position and velocity. Such events are negligible in a continuous phase-space setting, as they represent higher-order interaction events that are expected to occur with much lower probability than binary encounters in a dilute population regime. Moreover, this formulation would make the activation process depend on a highly specific alignment of microscopic velocities, whereas the ecological mechanisms motivating Case II are intended to reflect environmental perception and prey abundance at a coarser level.

The hybrid term $\alpha \rho_f f g_0$ can instead be interpreted as a mean-field closure of a nonlocal perception process, in which predators respond to environmental cues (e.g., visual density, chemical signals, or accumulated experience), rather than to independent microscopic contacts alone. More precisely, we interpret $\rho_f(x,t)$ as (proportional to) a coarse-grained environmental cue perceived by predators at an intermediate spatial scale, larger than the microscopic encounter scale but still much smaller than the macroscopic scale emerging in the hydrodynamic limit discussed in Section~\ref{sec:3}. In this sense, $\rho_f(x,t)$ acts as a local environmental field that modulates the success probability of activation upon contact. 

Within this interpretation, the local interaction term $f g_0$ continues to describe the physical encounter required for predator activation, whereas the density factor $\rho_f$ modulates the activation rate according to the perceived abundance of prey in the surrounding environment.

Equivalently, the hybrid interaction may be viewed as the local approximation of a nonlocal perceptive activation mechanism $\alpha f g_0(K * \rho_f)$, where $K$ describes a predator sensing kernel and $(K * \rho_f)(x) = \int\!\rmd y\, K(x-y)\rho_f(y)$ denotes the corresponding perception field. Accordingly, predators perceive prey abundance over a spatial range that is larger than the microscopic encounter scale but remains small compared with the macroscopic scale governing the hydrodynamic description. Under the localization of the sensing kernel, one formally recovers the local interaction term $\alpha \rho_f f g_0$ as the leading-order approximation of the perception mechanism. The local and nonlocal descriptions therefore yield the same hydrodynamic closure at the orders considered in the present work.

The linear dependence on $\rho_f$ should be regarded as the simplest perception law. More general perception functions could also be incorporated within the same framework, leading to different activation nonlinearities. An advantage of this choice, however, is that it already produces a sigmoidal functional response through the activation dynamics and the associated fast-state reduction, without requiring more elaborate perception mechanisms.

This modeling choice preserves the two-body interaction structure required for the physical occurrence of predation events, while introducing a minimal and biologically meaningful nonlocality through the macroscopic density. From a modeling perspective, such hybrid constructions naturally arise whenever individual behavior is assumed to depend on aggregate environmental quantities, and can be formally justified as limits of nonlocal interaction kernels or perception processes with finite sensing radius. In this sense, it should not be viewed as an ad hoc modification, but a way to enforce consistency between the microscopic interaction structure and the macroscopic behavioral mechanisms that give rise to sigmoidal functional responses.
\end{remark}

\section{Hydrodynamic limit}
\label{sec:3}

The kinetic formulation~\eqref{eq:prey} and \eqref{eq:pred1}--\eqref{eq:pred2} links movement, behavior, and predation at the individual scale. In this section, we derive hydrodynamic equations for the macroscopic closed system involving solely the ``slow variables'': the prey density $\rho_f$ and the total predator density $\rho_P := \rho_{g_0} + \rho_{g_1} + \rho_{g_2}$. This reduction is obtained by exploiting the fast relaxation in velocity via a diffusive rescaling of spacetime variables, combined with a suitable scaling of the internal behavioral transitions. This yields a reaction--diffusion system of Rosenzweig--MacArthur type, featuring reactive terms that produce Holling type II or III functional responses, equations \eqref{eq:RMII} and \eqref{eq:RMIII} below, depending on the specific activation mechanism considered.

We introduce the scaling parameter $0<\eps \ll 1$ and, by an abuse of notation, retain the same letters $(x,v,t)$ for the macroscopic variables (i.e., through the rescaling $(x,v,t)$ to $(\eps^{-1} x, v, \eps^{-2}t)$). In these new variables, the Liouville operator on the left-hand side of the kinetic equations becomes $\eps^2\partial_t + \eps v \cdot \nabla_x$. Due to the linear restoring force $-v$ in $\mc L_j$ ($j=1,2$) acting at high velocities, the velocity distribution undergoes a drift toward the origin, which counteracts the unbounded spreading induced by the Laplacian. Consequently, the Ornstein-Uhlenbeck dynamics admits the Maxwellian
\begin{equation}
\label{eq:MD}
M_j(v) = \frac1{(2\pi D_j)^{d/2}} \exp\left(-\frac{|v|^2}{2D_j}\right)
\end{equation}
as its unique equilibrium probability density, ensuring fast relaxation toward it on the $O(\varepsilon^2)$ timescale. The presence of such a restoring force is essential for the kinetic description to be consistently reduced to macroscopic equations for the local densities
\[
\rho_{f}(x,t) = \int\! \rmd v\,  f(x,v,t), \qquad \rho_{g_i}(x,t)=\int\! \rmd v\, g_i(x,v,t), \quad i=0,1,2,
\]
in the limit $\eps\to 0$.

The remaining terms on the right-hand side of \eqref{eq:prey} and \eqref{eq:pred1}--\eqref{eq:pred2} represent phenomena (activation, handling, prey logistic dynamics, and basal mortality) that occur on a slower timescale than velocity relaxation. This separation is implemented by choosing a suitable scaling with respect to the parameter $\eps$ for all the parameters involved. To this purpose, we observe that the internal transitions within the predator subsystem, $g_0 \rightarrow g_1 \rightarrow g_2 \rightarrow g_0$, occur on a fast timescale relative to the ecological dynamics, so that the internal sub-populations are expected to rapidly reach a quasi-stationary behavioral equilibrium. This is achieved by scaling the transition rates $(\beta, h^{-1})$ associated with hunting and handling as $(\eps\beta, \eps h^{-1})$. Crucially, the activation rate itself will take a linear form $\eps\alpha$ in Case I, and a density-dependent form $\eps\alpha\rho_f$ in Case II. Concurrently, the purely ecological parameters $(r,a,m, e a)$ are scaled as $(\eps^2r,\eps^2 a,\eps^2 m, \eps^2 e a)$ to match the slow macroscopic time horizon. Since the internal transitions of the predator subsystem preserve mass, the total predator distribution function $g_P := g_0 + g_1 + g_2$ evolves according to the equation
\[
\eps^2 \partial_t g_P + \eps v\cdot\nabla_x g_P= \mc L_2  g_P + \eps^2 e a f g_1 - \eps^2 m g_P.
\]
Therefore, in both cases, it is possible to describe the macroscopic dynamics via a closed system of equations for $\rho_f(x,t)$ and $\rho_P(x,t) = \int\!\rmd v\, g_P(x,v,t)$, where the specific functional response emerges naturally from the underlying behavioral structure. This is the content of the following subsections.

\subsection{Case I: Opportunistic activation (Holling type II)}
\label{sec:3.1}

In macroscopic variables, the system of equations \eqref{eq:prey}--\eqref{eq:pred1} reads:
\begin{equation}
\label{eq:prey-predII}
\begin{aligned}
\partial_t f + \frac 1\eps v\cdot\nabla_x f & = \frac1{\eps^2}\mc L_1 f + r f \Big(1-\frac{\rho_f}{K}\Big) - a f g_1, \\ 
\partial_t g_0 + \frac 1\eps v\cdot\nabla_x g_0 & = \frac1{\eps^2}\mc L_2  g_0 + \frac 1\eps \left(- \alpha f g_0 + h^{-1} g_2\right) - m g_0 + e a f g_1, \\ \partial_t g_1 + \frac 1\eps  v\cdot\nabla_x g_1 & = \frac1{\eps^2}\mc L_2  g_1 + \frac 1\eps \left(\alpha f g_0 - \beta f g_1\right) - m g_1, \\ \partial_t g_2 + \frac 1\eps v\cdot\nabla_x g_2 & = \frac1{\eps^2}\mc L_2  g_2 + \frac 1\eps \left( \beta f g_1 - h^{-1} g_2 \right) - m g_2.
\end{aligned}
\end{equation}
We apply the Hilbert method to derive the hydrodynamic equations. We expand the unknown distribution functions $(f, g_0, g_1, g_2)$---whose explicit dependence on the scaling parameter $\eps$ is hereafter omitted---in powers of $\eps$:
\[
F = F^{(0)} + \eps F^{(1)} + \eps^2 F^{(2)} + \dots, \quad \text{for } F \in \{f, g_0, g_1, g_2\},
\]
where the coefficients $F^{(k)}$ are independent of $\eps$. Substituting these expansions into the kinetic equations and identifying terms of the same order, we obtain a hierarchy of problems.
\begin{remark}
\label{rem:solvability}
The hierarchy equations take the form $\mc L_j \mc X = \mc B$, $j=1,2$, where $\mc B$ is a known term (possibly depending on $(x,t) \in \bb R^d\times \bb R_+$, which here play the role of parameters). In order for an equation of this form to admit a solution, $\mc B$ must satisfy the Fredholm solvability condition. To determine it, we recall that the adjoint Ornstein-Uhlenbeck operator satisfies $\ker(\mc L_j^*) = \text{span}\{1\}$, while the corresponding Fokker-Planck operator satisfies $\ker(\mc L_j) = \text{span}\{M_j\}$, with $M_j$ as in \eqref{eq:MD}, see, e.g., \cite{Villani2009}. Hence, $\mc B$ must be orthogonal to constants, namely $\int\!\rmd v\, \mc B = 0$. Under this condition, since $\ker(\mc L_j) = \text{span}\{M_j\}$, the equation admits infinitely many solutions, which differ by a multiple of $M_j$. We denote by $\mc L_j^{-1} \mc B$ the unique solution satisfying $\int\! \rmd v\, \mc L_j^{-1} \mc B = 0$. Therefore, the general solution reads $\mc X = \mc L_j^{-1} \mc B + \rho M_j$, where $\rho = \int\!\rmd v\, \mc X$.
\end{remark}
\noindent $\bullet$ Order $\eps^{-2}$: At the leading order, we find $\mc L_1 f^{(0)} = 0$ and $\mc L_2 g_i^{(0)} = 0$ for $i=0,1,2$, which imply
\begin{equation}
\label{eq:order0}
f^{(0)}(x,v,t) = \rho_f^{(0)}(x,t) M_1(v), \quad g_i^{(0)}(x,v,t) = \rho_{g_i}^{(0)}(x,t) M_2(v), \quad i=0,1,2,
\end{equation}
where we denote by $\rho_f^{(j)}(x,t)$ and $\rho_{g_i}^{(j)}(x,t)$ (and later $\rho_P^{(j)}$) the $j$-th order terms in the $\eps$-expansion of the local densities.
\\ $\bullet$ Order $\eps^{-1}$: At this order, we obtain the following system:
\begin{equation}
\label{eq:order1}
\begin{aligned}
v\cdot \nabla_x f^{(0)} & = \mc L_1 f^{(1)}, \\
v\cdot \nabla_x g_0^{(0)} & = \mc L_2 g_0^{(1)} + \left(- \alpha f^{(0)} g_0^{(0)} + h^{-1} g_2^{(0)} \right),  \\ v\cdot \nabla_x g_1^{(0)} & = \mc L_2 g_1^{(1)} +  \left(\alpha f^{(0)} g_0^{(0)} - \beta f^{(0)} g_1^{(0)}\right), \\v\cdot \nabla_x g_2^{(0)} & = \mc L_2 g_2^{(1)} + \left( \beta f^{(0)} g_1^{(0)} - h^{-1} g_2^{(0)} \right).
\end{aligned}
\end{equation}
As, by \eqref{eq:order0},
\[
\int\!\rmd v \, v\cdot \nabla_x f^{(0)} = \nabla_x  \rho_f^{(0)}(x,t) \cdot \int\!\rmd v \, v  M_1(v) = 0,
\] 
the first equation in \eqref{eq:order1} admits solution, and precisely\footnote{We use here and later on that, given $c\in \bb R^d$, $\mc L_j^{-1}(M_j v\cdot c) = - M_j v \cdot c$, as can be immediately verified by noticing that $\mc L_j f = D_j \nabla_v \cdot (M_j\nabla_v(f/M_j))$.} 
\begin{equation}
\label{eq:f1}
f^{(1)}(x,v,t) = - M_1(v) v \cdot \nabla_x\rho_f^{(0)}(x,t) + \rho_f^{(1)}(x,t) M_1(v).
\end{equation}
For the remaining equations, again by \eqref{eq:order0} we have
\[
\int\!\rmd v \, v\cdot \nabla_x g_i^{(0)} = \nabla_x  \rho_{g_i}^{(0)}(x,t) \cdot \int\!\rmd v \, v  M_2(v) = 0,
\] 
so that the solvability conditions reduce to 
\[
\int\!\rmd v\,  \left(- \alpha f^{(0)} g_0^{(0)} + h^{-1} g_2^{(0)} \right) = \int\!\rmd v\, \left(\alpha f^{(0)} g_0^{(0)} - \beta f^{(0)} g_1^{(0)}\right)   = \int\!\rmd v\, \left( \beta f^{(0)} g_1^{(0)} - h^{-1} g_2^{(0)} \right) = 0,
\]
which, letting $\kappa = \int\!\rmd v \, M_1(v) M_2(v) = [2\pi(D_1+D_2)]^{-d/2}$, leads to
\[
- \alpha \kappa \rho_f^{(0)} \rho_{g_0}^{(0)} + h^{-1} \rho_{g_2}^{(0)} = \alpha \kappa \rho_f^{(0)} \rho_{g_0}^{(0)} - \beta \kappa  \rho_f^{(0)} \rho_{g_1}^{(0)} = \beta \kappa  \rho_f^{(0)} \rho_{g_1}^{(0)} - h^{-1} \rho_{g_2}^{(0)} = 0.
\]
Assuming non-vanishing prey density ($\rho_f^{(0)}>0$), these algebraic relations imply that, at the leading order, the macroscopic densities of the predator states equilibrate. Specifically, they define a one-dimensional manifold in the $(\rho_{g_0}, \rho_{g_1}, \rho_{g_2})$ space, allowing the individual densities $\rho_{g_i}^{(0)}$ to be uniquely determined as functions of the total predator density $\rho_P^{(0)}$ and the prey density $\rho_f^{(0)}$. More precisely, solving these relations yields:
\begin{equation}
\label{algeb-equil}
\rho_{g_0}^{(0)} = \frac{\rho_P^{(0)}}{1 + \frac{\alpha}{\beta}+h\alpha\kappa\rho_f^{(0)}}, \qquad \rho_{g_1}^{(0)} = \frac{\frac{\alpha}{\beta}\rho_P^{(0)}}{1+\frac{\alpha}{\beta}+h\alpha\kappa\rho_f^{(0)}\,}, \qquad \rho_{g_2}^{(0)} =  \frac{h \alpha\kappa \rho_f^{(0)}\rho_P^{(0)}}{1+\frac{\alpha}{\beta} + h\alpha\kappa\rho_f^{(0)}}.
\end{equation}
Finally, the first-order corrections for the distribution functions are given by
\begin{equation}
\label{eq:g1}
g_i^{(1)}(x,v,t) = - M_2(v) v \cdot \nabla_x\rho_{g_i}^{(0)}(x,t) +  \rho_{g_i}^{(1)}(x,t) M_2(v) + \mc L_2^{-1} \Theta_i(x,v,t), \quad i=0,1,2,
\end{equation}
where $\Theta_i$ denotes the respective terms in parenthesis on the right-hand side of \eqref{eq:order1}.
\\ $\bullet$ Order $\eps^0$: At this order, the hierarchy reads:
\begin{equation}
\label{eq:order2}
\begin{aligned}
\partial_t f^{(0)} + v\cdot\nabla_x f^{(1)} & = \mc L_1 f^{(2)} + r f^{(0)} \Big(1-\frac{\rho_f^{(0)}}{K}\Big) - a f^{(0)} g_1^{(0)}, \\ 
\partial_t g_0^{(0)} + v\cdot\nabla_x g_0^{(1)} & = \mc L_2  g_0^{(2)} + \left[- \alpha \left(f^{(0)} g_0^{(1)} + f^{(1)} g_0^{(0)}\right) + h^{-1} g_2^{(1)}\right] - m g_0^{(0)} + e a f^{(0)} g_1^{(0)}, \\ \partial_t g_1^{(0)} + v\cdot\nabla_x g_1^{(1)} & = \mc L_2  g_1^{(2)} + \left[\alpha \left(f^{(0)} g_0^{(1)} + f^{(1)} g_0^{(0)}\right) - \beta\left(f^{(0)} g_1^{(1)} + f^{(1)} g_1^{(0)}\right)\right] - m g_1^{(0)}, \\ \partial_t g_2^{(0)} + v\cdot\nabla_x g_2^{(1)} & = \mc L_2  g_2^{(2)} + \left[\beta\left(f^{(0)} g_1^{(1)} + f^{(1)} g_1^{(0)}\right)\ - h^{-1} g_2^{(1)} \right] - m g_2^{(0)}.
\end{aligned}
\end{equation}
 For the first equation to admit a solution $f^{(2)}$, the solvability condition requires
\[
\partial_t\rho_f^{(0)} + \int\!\rmd v\, v\cdot\nabla_x f^{(1)} =  r \rho_f^{(0)} \Big(1-\frac{\rho_f^{(0)}}{K}\Big) - a\kappa \rho_f^{(0)} \rho_{g_1}^{(0)}.
\]
Recalling \eqref{eq:f1}, the integral term yields the diffusion flux\footnote{Note that the component of $f^{(1)}$ belonging to $\ker(\mc L_1)$ does not contribute to the flux since $\int\!\rmd v\, v M_1(v) = 0$.},
\[
\int\!\rmd v\, v\cdot\nabla_x f^{(1)} = - \int\!\rmd v\, M_1(v)(v\cdot\nabla_x)  (v \cdot \nabla_x\rho_f^{(0)}) = - \int\!\rmd v\, M_1(v) \mathrm{trace} [v \otimes v \, \partial^2\rho_f^{(0)}]  = - D_1 \Delta \rho_f^{(0)}.
\]
Analogous conditions are required for the existence of $g_i^{(2)}$. While we omit their explicit form for each $i$, we focus on the evolution of the total predator density $\rho_P^{(0)}$. Summing the three predator equations and noting that the sum of the three terms in square brackets on the right-hand side of \eqref{eq:order2} is zero, we obtain 
\[
\partial_t \rho_P^{(0)} +  \int\!\rmd v\, v\cdot\nabla_x g_P^{(1)} = ea\kappa \rho_f^{(0)} \rho_{g_1}^{(0)} - m \rho_P^{(0)}.
\]
Since the sum of the $\Theta_i$'s appearing in the right-hand side of \eqref{eq:g1} also vanishes, the total first-order distribution is simply $g_P^{(1)}(x,v,t) = - M_2(v) v \cdot \nabla_x\rho_P^{(0)}(x,t) + \rho_P^{(1)}(x,t) M_2(v)$, which implies
\[
\int\!\rmd v\, v\cdot\nabla_x g_P^{(1)} = - D_2 \Delta \rho_P^{(0)}.
\]

In conclusion, in the limit $\eps\to 0$, the time evolution of the slow variables $(\rho_f,\rho_P) = (\rho_f^{(0)},\rho_P^{(0)})$ is governed by the following  reaction diffusion system:
\begin{equation}
\label{eq:RMII}
\begin{cases}
\partial_t \rho_f = D_1 \Delta \rho_f + {\displaystyle r \rho_f\Big(1 - \frac{\rho_f}{K}\Big)}  - \Phi_{\rm II}(\rho_f,\rho_P), \\ \partial_t \rho_P = D_2 \Delta \rho_P + e \Phi_{\rm II}(\rho_f,\rho_P) - m  \rho_P,
\end{cases}
\end{equation}
where, in view of the algebraic relations \eqref{algeb-equil}, the predation flux $\Phi(\rho_f,\rho_P) = a\kappa \rho_f\rho_{g_1}^{(0)}$ takes the form of a  Holling type~II functional response: 
\[
\Phi(\rho_f,\rho_P) = \Phi_{\rm II}(\rho_f,\rho_P) = \frac{A_2\rho_f \rho_P}{1+H_2\rho_f}, \qquad \text{ with }
A_2 = \frac{\alpha\kappa a}{\alpha+\beta} \text{ and } H_2 =\frac{\alpha\beta h\kappa}{\alpha+\beta}.
\]

\begin{remark}
\label{rem:3.1}
It is worth noting that the macroscopic system \eqref{eq:RMII} describes the dynamics only after fast transient processes associated with the diffusive scaling. Indeed, since the initial conditions for the kinetic distributions do not generally satisfy the equilibrium constraints \eqref{eq:order0} and \eqref{algeb-equil}, an initial relaxation layer occurs before the asymptotic regime is reached. During this stage, the velocity distributions relax toward local Maxwellians and the predator sub-populations approach the equilibrium manifold determined by the prey and total predator densities. A treatment of these transients would require a separate multi-layer analysis using stretched time variables in order to match the initial data with the Hilbert expansion.
\end{remark}

\subsection{\bf Case II: Learning-driven activation (Holling type III)}
\label{sec:3.2}

In macroscopic variables, the system of equations \eqref{eq:prey}--\eqref{eq:pred2} takes the form:
\begin{equation}
\label{eq:prey-predIII}
\begin{aligned}
\partial_t f + \frac 1\eps v\cdot\nabla_x f & = \frac1{\eps^2}\mc L_1 f + r f \Big(1-\frac{\rho_f}{K}\Big) - a f g_1, \\ 
\partial_t g_0 + \frac 1\eps v\cdot\nabla_x g_0 & = \frac1{\eps^2}\mc L_2  g_0 + \frac 1\eps \left(- \alpha \rho_f f g_0 + h^{-1} g_2\right) - m g_0 + e a f g_1, \\ \partial_t g_1 + \frac 1\eps  v\cdot\nabla_x g_1 & = \frac1{\eps^2}\mc L_2  g_1 + \frac 1\eps \left(\alpha \rho_f  f g_0 - \beta f g_1\right) - m g_1, \\ \partial_t g_2 + \frac 1\eps v\cdot\nabla_x g_2 & = \frac1{\eps^2}\mc L_2  g_2 + \frac 1\eps \left( \beta f g_1 - h^{-1} g_2 \right) - m g_2.
\end{aligned}
\end{equation}
We proceed as before, by applying the Hilbert method and obtaining a hierarchy of problems:
\\ $\bullet$ Order $\eps^{-2}$: At the leading order, we find $\mc L_1 f^{(0)} = 0$ and $\mc L_2 g_i^{(0)} = 0$ for $i=0,1,2$, which imply \eqref{eq:order0}.
\\ $\bullet$ Order $\eps^{-1}$: At this order, we obtain the following system:
\begin{equation}
\label{eq:order3}
\begin{aligned}
v\cdot \nabla_x f^{(0)} & = \mc L_1 f^{(1)}, \\
v\cdot \nabla_x g_0^{(0)} & = \mc L_2 g_0^{(1)} + \left(- \alpha \rho_f^{(0)} f^{(0)} g_0^{(0)} + h^{-1} g_2^{(0)} \right),  \\ v\cdot \nabla_x g_1^{(0)} & = \mc L_2 g_1^{(1)} +  \left(\alpha \rho_f^{(0)} f^{(0)} g_0^{(0)} - \beta f^{(0)} g_1^{(0)}\right), \\v\cdot \nabla_x g_2^{(0)} & = \mc L_2 g_2^{(1)} + \left( \beta f^{(0)} g_1^{(0)} - h^{-1} g_2^{(0)} \right).
\end{aligned}
\end{equation}
By arguing as in the previous case, the first equation has solution \eqref{eq:f1}, while for the remaining equations to admit solutions $g_i^{(1)}$, the terms on the right-hand side must satisfy the orthogonality conditions
\[
\begin{aligned}
& \int\!\rmd v\,  \left(- \alpha \rho_f^{(0)} f^{(0)} g_0^{(0)} + h^{-1} g_2^{(0)} \right) = \int\!\rmd v\, \left(\alpha \rho_f^{(0)} f^{(0)} g_0^{(0)} - \beta f^{(0)} g_1^{(0)}\right)  \\ &  \qquad = \int\!\rmd v\, \left( \beta f^{(0)} g_1^{(0)} - h^{-1} g_2^{(0)} \right) = 0,
\end{aligned}
\]
which lead to
\[
- \alpha \kappa (\rho_f^{(0)})^2 \rho_{g_0}^{(0)} + h^{-1} \rho_{g_2}^{(0)} = \alpha \kappa (\rho_f^{(0)})^2 \rho_{g_0}^{(0)} - \beta \kappa  \rho_f^{(0)} \rho_{g_1}^{(0)} = \beta \kappa  \rho_f^{(0)} \rho_{g_1}^{(0)} - h^{-1} \rho_{g_2}^{(0)} = 0.
\]
This algebraic system defines the quasi-stationary behavioral equilibrium for the predator sub-popula\-tions at the leading order. Assuming a non-vanishing prey density ($\rho_f^{(0)} > 0$),  solving these relations yields
\begin{equation}
\label{algeb-equilbis}
\begin{aligned}
\rho_{g_0}^{(0)} & = \frac{\rho_P^{(0)}}{1 + \frac{\alpha}{\beta}\rho_f^{(0)} + h\alpha\kappa(\rho_f^{(0)})^2}, \qquad \rho_{g_1}^{(0)} = \frac{\frac{\alpha}{\beta}\rho_f^{(0)}\rho_P^{(0)}}{1 + \frac{\alpha}{\beta}\rho_f^{(0)} + h\alpha\kappa(\rho_f^{(0)})^2}, \\ \rho_{g_2}^{(0)} & =  \frac{h \alpha\kappa(\rho_f^{(0)})^2\rho_P^{(0)}}{1+\frac{\alpha}{\beta}\rho_f^{(0)} + h\alpha\kappa(\rho_f^{(0)})^2}.
\end{aligned}
\end{equation}
Finally, the first-order corrections for the distribution functions are given by:
\begin{equation}
\label{eq:g1bis}
g_i^{(1)}(x,v,t) = - M_2(v) v \cdot \nabla_x\rho_{g_i}^{(0)}(x,t) +  \rho_{g_i}^{(1)}(x,t) M_2(v) + \mc L_2^{-1} \bar\Theta_i(x,v,t), \quad i=0,1,2,
\end{equation}
where $\bar\Theta_i$ denotes the respective terms in parenthesis on the right-hand side of \eqref{eq:order3}.
\\ $\bullet$ Order $\eps^{0}$: By repeating the same arguments used in Case I, since also in this case the fast transition terms telescope to zero, the equations for the prey and total predator distributions read:
\begin{equation}
\label{eq:order2_correct}
\begin{aligned}
\partial_t f^{(0)} + v\cdot\nabla_x f^{(1)} & = \mc L_1 f^{(2)} + r f^{(0)} \Big(1-\frac{\rho_f^{(0)}}{K}\Big) - a f^{(0)} g_1^{(0)}, \\
\partial_t g_P^{(0)} + v\cdot\nabla_x g_P^{(1)} & = \mc L_2 g_P^{(2)} + e a f^{(0)} g_1^{(0)} - m g_P^{(0)}.
\end{aligned}
\end{equation}
We project the kinetic equations \eqref{eq:order2_correct} onto $\text{ker}(\mc L_1^*) = \text{span}\{1\}$ and $\text{ker}(\mc L_2^*) = \text{span}\{1\}$ by integrating over the velocity space $\mathbb{R}^d$. We thus obtain
\[
\begin{aligned}
\partial_t \rho_f^{(0)} - D_1 \Delta \rho_f^{(0)} & = r \rho_f^{(0)} \left(1 - \frac{\rho_f^{(0)}}{K}\right) - a \kappa \rho_f^{(0)} \rho_{g_1}^{(0)}, \\ 
\partial_t \rho_P^{(0)} - D_2 \Delta \rho_P^{(0)} & = e a \kappa \rho_f^{(0)} \rho_{g_1}^{(0)} - m \rho_P^{(0)}.
\end{aligned}
\]
To derive the above equations, we argued as in the Case I, by using \eqref{eq:f1} and \eqref{eq:g1bis}  (where the sum of the $\bar \Theta_i$'s vanishes) to compute the diffusion fluxes.

In conclusion, in the limit $\eps\to 0$, the time evolution of the slow variables $(\rho_f,\rho_P) = (\rho_f^{(0)},\rho_P^{(0)})$ is governed by the reaction diffusion system:
\begin{equation}
\label{eq:RMIII}
\begin{cases}
\partial_t \rho_f = D_1 \Delta \rho_f + {\displaystyle r  \rho_f\Big(1 - \frac{\rho_f}{K}\Big)}  - \Phi_{\rm III}(\rho_f,\rho_P), \\ \partial_t \rho_P = D_2 \Delta \rho_P + e \Phi_{\rm III}(\rho_f,\rho_P) - m  \rho_P,
\end{cases}
\end{equation}
where, in view of the algebraic relation \eqref{algeb-equilbis}, the predation flux $\Phi(\rho_f,\rho_P) = a\kappa \rho_f\rho_{g_1}$ takes here the form
\[
\Phi(\rho_f, \rho_P) = \Phi_{\rm III}(\rho_f, \rho_P) = \frac{\alpha \kappa a \rho_f^2 \rho_P}{\beta + \alpha \rho_f + \alpha \beta h \kappa \rho_f^2},
\]
which represents a generalized Holling type~III functional response. Matching both the quadratic onset $\Phi\sim (a\alpha\kappa/\beta)\rho_f^2\rho_P$ as $\rho_f\to 0$ and the saturating level $\Phi\to a (\beta h)^{-1} \rho_P$ as $\rho_f\to\infty$ with the canonical Holling type~III form, neglecting the intermediate linear density effect $\alpha \rho_f$ in the denominator, we can also approximate
\begin{equation}
\label{eq:Phi-typeIII} 
\Phi_{\rm III}(\rho_f,\rho_P) \approx \frac{A_3 \rho_f^2 \rho_P}{1 + H_3\rho_f^2}  \qquad \text{ with }
A_3 = \frac{\alpha\kappa a}{\beta} \text{ and } H_3 = \alpha \kappa h.
\end{equation}

\section{Discussion and concluding remarks}
\label{sec:4}

A key qualitative difference between the functional responses derived in Sections~\ref{sec:3.1}--\ref{sec:3.2} lies in their behavior at low prey density. This property, while classical at the macroscopic level, acquires a direct mechanistic interpretation within the present kinetic framework.

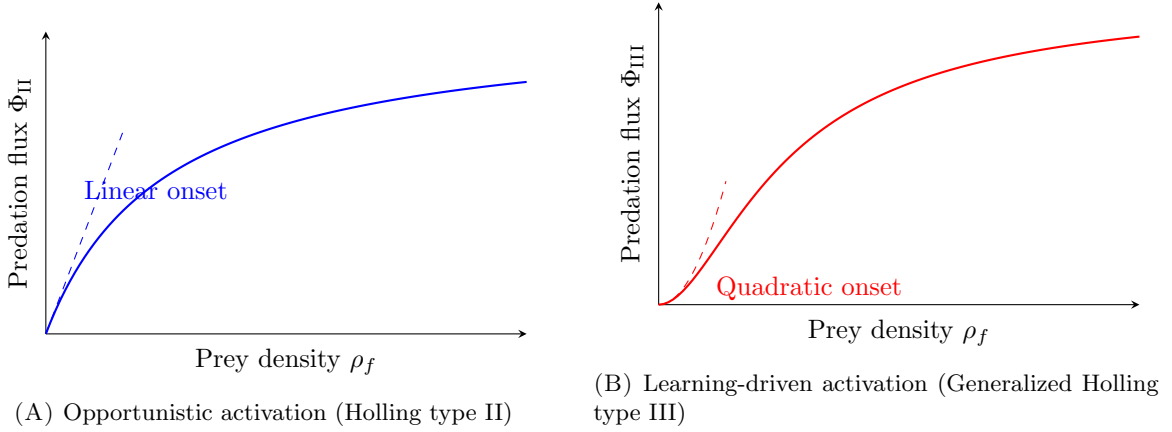
\begin{figure}[htbp]
    \centering
    \begin{subfigure}[b]{0.48\textwidth}
        \centering
        \begin{tikzpicture}
        \begin{axis}[
            axis lines = left,
            xlabel = {Prey density $\rho_f$},
            ylabel = {Predation flux $\Phi_{\rm II}$},
            xmin=0, xmax=5,
            ymin=0, ymax=1.2,
            xtick=\empty, ytick=\empty,
            width=0.85\textwidth,
            height=4cm,
            scale only axis, 
            clip=false
        ]
            \addplot [
                domain=0:5, 
                samples=100, 
                color=blue, 
                thick
            ] {x / (1 + 0.8*x)};
            
            \addplot [
                domain=0:0.8, 
                samples=2, 
                color=blue, 
                dashed
            ] {x};
            
            \node[color=blue, anchor=south west] at (axis cs: 0.3,0.5) {Linear onset};
        \end{axis}
        \end{tikzpicture}
        \caption{Opportunistic activation (Holling type II)}
        \label{fig:holling2}
    \end{subfigure}
    \hfill
    \begin{subfigure}[b]{0.48\textwidth}
        \centering
        \begin{tikzpicture}
        \begin{axis}[
            axis lines = left,
            xlabel = {Prey density $\rho_f$},
            ylabel = {Predation flux $\Phi_{\rm III}$},
            xmin=0, xmax=5,
            ymin=0, ymax=1.2,
            xtick=\empty, ytick=\empty,
            width=0.85\textwidth,
            height=4cm,
            scale only axis, 
            clip=false
        ]
            \addplot [
                domain=0:5, 
                samples=100, 
                color=red, 
                thick
            ] {x^2 / (1 + 0.5*x + 0.8*x^2)};
            
            \addplot [
                domain=0:0.7, 
                samples=50, 
                color=red, 
                dashed
            ] {x^2};
            
            \node[color=red, anchor=north west] at (axis cs: 0.5,0.15) {Quadratic onset};
        \end{axis}
        \end{tikzpicture}
        \caption{Learning-driven activation (Generalized Holling type III)}
        \label{fig:holling3}
    \end{subfigure}
    
    \caption{Qualitative comparison of the derived functional responses at low prey density $\rho_f \to 0$. In the opportunistic scenario (A), the flux exhibits a linear tangent at the origin. Conversely, the learning-driven mechanism (B) produces a sigmoidal profile with a quadratic onset at low prey density.}
    \label{fig:functional_responses}
\end{figure}

In the opportunistic activation scenario (Case~I), the functional response satisfies
\[
\Phi_{\rm II}(\rho_f,\rho_P) \sim A_2 \rho_f \rho_P \quad \text{as } \rho_f \to 0,
\]
indicating that predation occurs at leading order even for arbitrarily small prey densities, see Figure~\ref{fig:functional_responses}(\subref{fig:holling2}). This linear regime reflects the fact that activation is driven by pairwise encounters, and thus a single prey is sufficient to trigger predator activity.

In contrast, in the learning-driven activation scenario (Case~II), one has
\[
\Phi_{\rm III}(\rho_f,\rho_P) \sim \frac{a\alpha\kappa}{\beta}\, \rho_f^2 \rho_P \quad \text{as } \rho_f \to 0,
\]
which implies a quadratic onset of the predation rate, see Figure~\ref{fig:functional_responses}(\subref{fig:holling3}). At the macroscopic level, this behavior may be interpreted as an effective threshold-like mechanism for predator establishment, since the predation term becomes negligible at low prey density.

This comparison shows that the qualitative differences between Holling type~II and type~III responses reflect distinct microscopic activation processes. In particular, the transition from a linear to a quadratic low-density regime corresponds to a shift from purely local interaction-driven activation to a mechanism modulated by environmental prey availability. As already noticed, this latter mechanism should not necessarily be interpreted as a specific cognitive process. Rather, it may be regarded as the result of a broad class of environmentally mediated activation mechanisms through which prey availability affects predator efficiency.

This study complements the work of Dawes and Souza \cite{DawesSouza2013}, who obtained Holling type II and III functional responses from stochastic individual-level models. By contrast, the present approach relies on a kinetic phase-space description incorporating spatial transport, velocity dynamics, and internal behavioral states, from which its macroscopic limit gives rise to the Rosenzweig--MacArthur-type reaction--diffusion systems. A distinctive feature of the learning-driven model is the hybrid micro-macro activation term discussed in Remark \ref{rem:2.2}, which provides a minimal representation of environmental modulation of predator activation while preserving the kinetic interaction structure.

We stress that the hydrodynamic limit presented here is formal. A rigorous justification would require combining the hypocoercive properties of the Ornstein-Uhlenbeck operators with compactness techniques for the nonlinear reaction terms and the fast internal-state transitions. While closely related diffusion limits have been established for several kinetic Fokker-Planck models, the present system appears to fall outside the scope of currently available theories. We leave this problem for future investigation.

Some of the underlying mechanisms related to the formal hydrodynamic limit developed in Section~\ref{sec:3}, can be explored numerically. More precisely, we consider numerical experiments for the kinetic and macroscopic systems in the one-dimensional setting on the periodic domain $[-L,L]$ with $L=10$.\footnote{Artificial intelligence tools were used to assist in the development and debugging of the Julia codes employed for these illustrative simulations.} The parameters are chosen as $D_1=D_2=0.1$, $r=K=a=1$, $m=0.2$, $e=0.5$, $\alpha=\beta=2$, $h=0.5$ and $\eps=0.05$. We remark that purpose of these simulations is not a quantitative validation of the limit, but rather an illustration of the fast relaxation towards the equilibrium manifold and the subsequent emergence of the diffusive regime.

To emphasize the differences between the opportunistic and the learning-driven activation mechanisms, we consider an initial prey density containing extended low-density regions, say
\[
\rho_f^0(x)=0.01+0.15\,\mathrm{e}^{-x^2},
\]
and an initial predator density
\[
\rho_P^0(x)=0.4+0.1\cos\!\left(\frac{\pi x}{L}\right).
\]

\begin{figure}[htbp] 
\centering 
\begin{subfigure}[b]{0.48\textwidth} \centering \includegraphics[width=\textwidth]{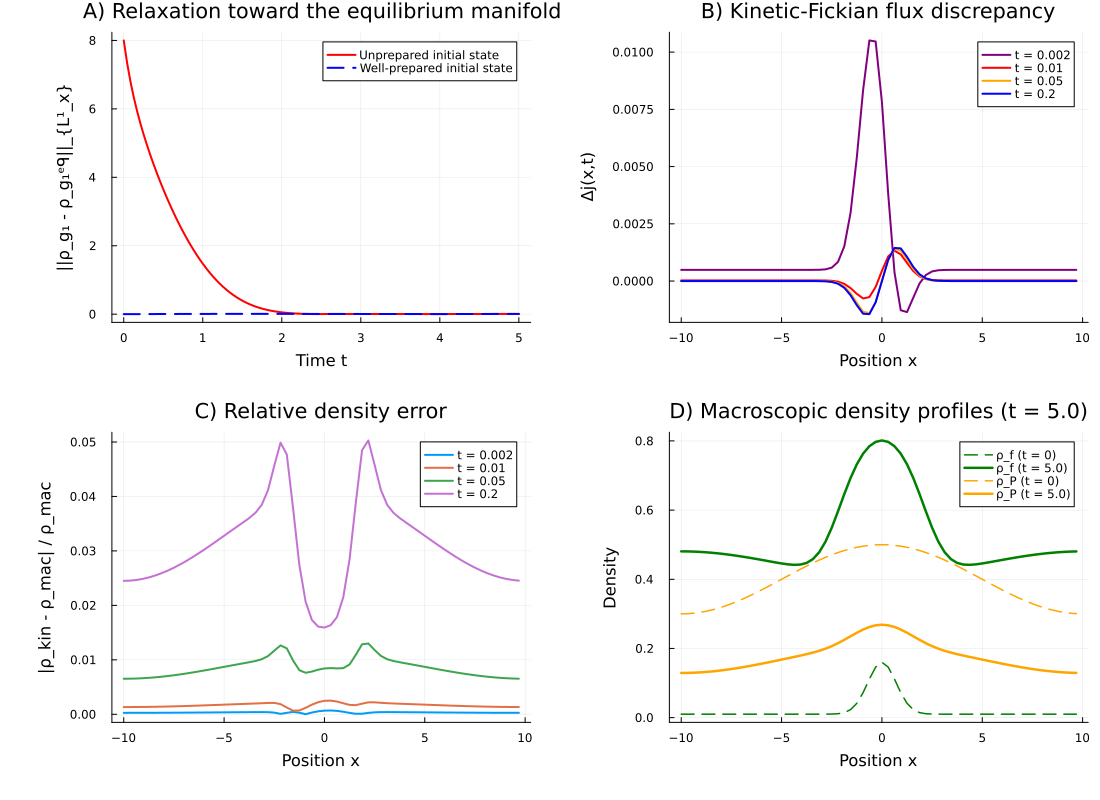} 
\caption{Case I: opportunistic activation.} \label{fig:num_caseI} 
\end{subfigure} \hfill 
\begin{subfigure}[b]{0.48\textwidth} \centering \includegraphics[width=\textwidth]{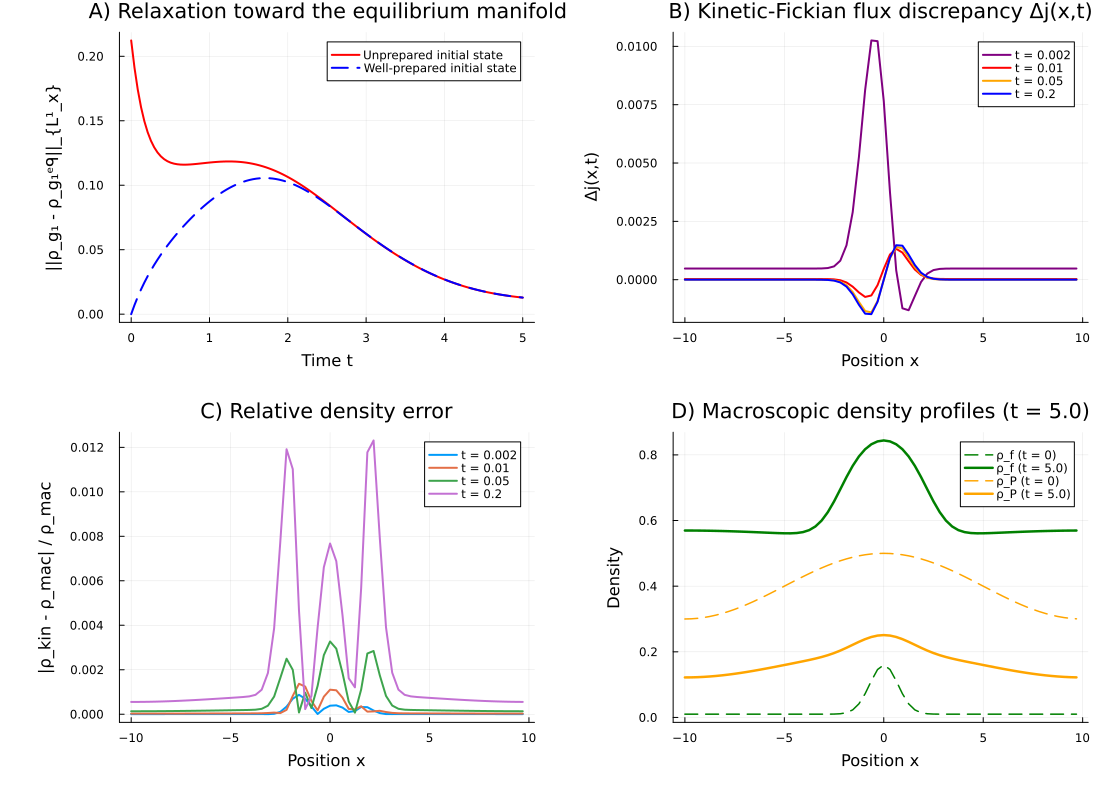} 
\caption{Case II: learning-driven activation.} \label{fig:num_caseII} 
\end{subfigure} 
\caption{Numerical illustration of the relaxation mechanisms underlying the formal hydrodynamic reduction for $\eps=0.05$. The simulations are initialized with $\rho_f^0(x) = 0.01+0.15\,\mathrm{e}^{-x^2}$ and $\rho_P^0(x) = 0.4+0.1 \cos (\pi x/L)$. Panel A displays the relaxation towards the equilibrium manifold, Panel B shows the flux discrepancy $\Delta j(x,t)=j_{\mathrm{kin}}(x,t)-j_{\mathrm{Fick}}(x,t)$, where $j_{\mathrm{kin}}(x,t)=\eps^{-1}\int\!\rmd v\, v\,f(x,v,t)$ and $j_{\mathrm{Fick}}(x,t)=-D_1\nabla_x\rho_f(x,t)$, Panel C reports the relative difference between kinetic and macroscopic prey densities, while Panel D compares the corresponding prey and predator density profiles. The simulations are intended as a qualitative illustration of the relaxation and diffusion mechanisms underlying the hydrodynamic reduction rather than as a quantitative convergence study.}
\label{fig:numerical_illustration} 
\end{figure}

We analyze two different kinetic distributions compatible with the prescribed macroscopic densities, which differ only in the preparation of the predator internal states. In the first one, all predators initially belong to the inactive class, whereas in the second one the predator internal states are initialized on the leading-order equilibrium manifold derived in Sections~\ref{sec:3.1} and \ref{sec:3.2}. In both configurations, the velocity distributions are set equal to Maxwellian profiles, while the prey distribution includes a small perturbation generating a non-zero initial flux. The simulations show the transient dynamics associated with the equilibration of the predator internal states and the subsequent emergence of the macroscopic regime. Since the velocity relaxation induced by the Fokker--Planck operator takes place within an initial layer of order $\eps^2$, the diagnostics displayed in Figure~\ref{fig:numerical_illustration} mainly capture the subsequent evolution of the associated macroscopic densities. In particular, Panel~A highlights the role of the internal-state preparation: initial data satisfying the leading-order quasi-equilibrium relations remain considerably closer to the equilibrium manifold throughout the evolution than the corresponding unprepared configurations. In the learning-driven case, the relaxation is less immediate, consistently with the explicit dependence of the quasi-equilibrium manifold on the prey density.

We conclude by remarking that several extensions of the model developed in this work can be considered. A first one consists in replacing the local activation factor $\rho_f$ by a nonlocal perception field $K * \rho_f$ associated with finite sensing radius. As discussed in Remark~\ref{rem:2.2}, the local formulation adopted in this work corresponds to the leading-order contribution of such a mechanism in the diffusive regime. Further extensions include the incorporation of additional behavioral states, state-dependent mobility coefficients, nonlocal interaction kernels, or stochastic perturbations.

\subsection*{Acknowledgements}

The author would like to thank Michela Ottobre for insightful discussions regarding the problem considered here.

\section*{Statements and declarations}

\begin{itemize}
\item \textbf{Funding} No funding was received for conducting this study.
\item \textbf{Conflict of interest} The author has no competing interests to declare that are relevant to the content of this article.
\item \textbf{Data Availability Statement} No datasets were generated or analysed during the current study. The Julia codes used to produce the numerical illustrations are available from the author upon reasonable request.
\end{itemize}


\begin{thebibliography}{99}

\bibitem{Beddington1975}
Beddington JD (1975) Mutual interference between parasites or predators and its effect on searching efficiency. J. Anim Ecol 44:331--340

\bibitem{Berardo2020} 
Berardo C, Geritz S, Gyllenberg M, Raoul G (2020) Interactions between different predator–prey states: a method for the derivation of the functional and numerical response. J Math Biol 80:2431--2468

\bibitem{BondesanZanella2025}
Bondesan A, Menale M, Toscani G, Zanella M (2025) Lotka--Volterra-type kinetic equations for interacting species. Nonlinearity 38, 075026 (35pp) 

\bibitem{ConfortoDesvillettesSoresina2017}
Conforto F, Desvillettes L, Soresina C (2018) About reaction--diffusion systems involving the Holling-type II and the Beddington-DeAngelis functional responses for predator--prey models. Nonlinear Differ Equ Appl 25, article 25

\bibitem{DawesSouza2013}
Dawes JHP, Souza MO (2013) A derivation of Holling’s type I, II and III functional responses in predator--prey systems. J Theor Biol 327(1):11-22

\bibitem{DeAngelis1975}
DeAngelis DL, Goldstein RA, Vandermeer R (1975) A model for trophic interaction. Ecology 56(4):881--892

\bibitem{DucrotLiuMagal2021}
Ducrot A, Liu Z, Magal P (2021) Large speed traveling waves for the Rosenzweig--MacArthur predator--prey model with spatial diffusion. Physica D 415:132730

\bibitem{Dunbar1986}
Dumbar RD (1986) Traveling waves in diffusive predator-prey equations: periodic orbits and point-to-periodic heteroclinic orbits. SIAM J Appl Math  46(6):1057--1078.

\bibitem{Holling1959}
Holling CS (1959) Some characteristics of simple types of predation and parasitism. Can Entomol 91(7):385--398

\bibitem{LiWu2008}
Li W, Wu S (2008) Traveling waves in a diffusive predator–prey model with holling type-III functional response. Chaos Solitons Fractals 37:476--486

\bibitem{Murray1989}
Murray JD (1989) Mathematical Biology. Springer--Verlag, New York--Berlin

\bibitem{Perthame2007}
Perthame B (2007) Transport Equations in Biology. Birkh\"auser, Basel

\bibitem{Real1977}
Real LA (1977) The kinetics of functional response. Am Nat 111(978):289--300

\bibitem{RosenzweigMacArthur1963}
Rosenzweig ML, MacArthur RH (1963) Graphical representation and stability conditions of predator--prey interactions. Am. Nat. 97:209--223

\bibitem{Smoller1994}
Smoller J (1994) Shock Waves and Reaction--Diffusion Equations. Springer-Verlag, New York, 2nd ed.

\bibitem{Villani2009}
Villani C (2009) Hypocoercivity. Memoirs of the American Mathematical Society 202(950)

\end{thebibliography}
\end{document}